\documentclass[preprint, 12pt]{elsarticle}
\usepackage{geometry} 
\usepackage{amssymb}
\usepackage{amsthm}
\usepackage{amsmath}
\usepackage{graphicx}
\usepackage{mathrsfs}
\usepackage{srcltx}
\newtheorem{theorem}{Theorem}
\newtheorem{definition}{Definition}
\newtheorem{problem}{Problem}

\begin{document}

		
		
		\title{The Bohr Radius of the Weighted Bloch Spaces}
		
		
		\author{Khasyanov R.Sh.\fnref{label2}}
		\ead{st070255@student.spbu.ru}
		\affiliation[label2]{organization={Saint Petersburg State University},
			addressline={University Avenue, 28D}, 
			city={Saint Petersburg},
			postcode={198504}, 
			country={Russian Federation}}
		
		\begin{abstract}
			The concept of the Bohr radius of a pair of Banach spaces is introduced. The lower estimate for the value of the Bohr radius from the Bloch space to the space of bounded functions  obtained by I. Kayumov, S. Ponnusamy and N. Shakirov is slightly improved. It is shown that for any weighted Bloch space the Bohr radius is not less than $1/\sqrt{2}$. A criterion for the sharpness of the Bohr inequality in the weighted Bloch space with $R=1 / \sqrt{2}$ is obtained. Using this criterion, the examples of the weights for which the inequality is sharp is given.
			
		\end{abstract}


	\maketitle
	
	\section{Introduction and main results}
	
	\vspace{1mm}
	In 1914, H. Bohr, studying Dirichlet series, noticed \cite{Bohr} the following interesting fact in complex analysis, which is now called the Bohr phenomenon:
	
	\newtheorem*{A}{Theorem A}\label{Th:A}
		\begin{A}
		Let $f(z)=\sum_{n\ge 0}a_nz^n$ and $\|f\|_{\infty}:=\sup_{z\in \mathbb{D}}|f(z)|\le1$ in the unit disc  $\mathbb{D}=\{|z|<1\}$. Then
		$$\sum_{n\ge0}^{} |a_n| r^n\le1,
		\quad 0\le r\le 1/3.$$
		The constant $1/3$ is sharp.
	\end{A}
	
	In fact, Bohr proved this for $r\le1/6$. The best constant $1/3$ was obtained independently in the same year
	by M.  Riesz, I. Schur and F. Wiener (for different proofs see e.g. appendix in \cite{Dixon05}).
	This theorem is equivalent to the inequality known as the Bohr inequality:
	$$
	\sum_{n\ge0}|a_n|r^n\le \|f\|_\infty,  \quad 0\le r\le 1/3.
	$$   
	
	An active study of various modifications and generalizations of the Bohr inequality began in the middle of 1990s since P. Dixon, using the Bohr inequality, solved a long-standing problem about the characterization of Banach algebras \cite{Dixon05}. Part of the subsequent research in this area is directed towards extending the Bohr phenomenon in multidimensional framework and in more abstract settings (see e.g. \cite{Aiz}, \cite{Boas}, \cite{BoasKhav}, \cite{PPS}). In 2018 B. Bhowmik and N. Das applied the Bohr inequality to the question of comparing majorant series of subordinating functions \cite{BhowDas}. For studies related to the Bohr inequality, see, for example, \cite{AHal}, \cite{BhowDas}-\cite{BDK}, \cite{GMR}-\cite{PPS}.
	
	\vspace{3mm}
	
	We consider the following questions:

	\begin{problem}\label{P:1}
		Fix a holomorphic function $f(z)=\sum_{n\ge0}a_nz^n$ and consider the class $\mathcal{F}$ of all analytic functions of the form  $\sum_{n\ge0}b_nz^n$ such that $|b_n|=|a_n|$ for all $n\ge 0.$ Which properties of $f$ are inherited by all functions in the class $\mathcal{F}$? In other words, what properties of a holomorphic function can be detected from the moduli of its Maclaurin series coefficients?
	\end{problem} 
	
	For instance, the radius of convergence of the Maclaurin series and the  Hardy  space $H^2$ norm depend only on the moduli of the coefficients. On the other hand,  there are examples of properties that may hold for some members of $\mathcal{F}$ but not for others: the series $\sum_{n\ge 1} \dfrac{\pm z^n}{n}$ are bounded in the $\mathbb{D}$ for almost every choice of plus and minus signs (\cite{Zyg}, Chapter V, Theorem 8.34), but if all plus signs are taken, the series is unbounded in the unit disk.
	Bohr's theorem implies that if $|f(z)|\le 1$ in the unit disc $\mathbb{D}$, then $|g(z)|\le 1$ in the disk of radius $1/3$ for all functions $g\in \mathcal{F}$.

	\begin{problem}\label{P:2}
		Let $X, Y$ be Banach spaces of analytic functions in the disc $\mathbb{D}=\left\{|z|<1\right\}$. Let $f(z)=\sum_{n\ge 0}a_nz^n$. We need to find the maximal $R$ for which
		$$
		\|f\|_{X}\le 1 \Longrightarrow  
		\|\sum_{n\ge 0}|a_n|(Rz)^n\|_{Y}\le1.
		$$ 
		We call such $R$ \textbf{the Bohr radius from $X$ to $Y$} and denote  $R_{X\rightarrow Y}$. If $X$ and $Y$ coincide, we write $R_X.$
	\end{problem}

	The classical Bohr theorem states that $R_{H^{\infty}}=1/3.$ We study Bohr radius in Bloch spaces.
	
	\begin{definition}\label{D:1}
		We say that the analytical function in $\mathbb{D}$ $f(z)=\sum_{n\ge0}a_nz^n$ belongs to the weighted Bloch space $\mathcal{B}(\omega)$ with the weight $\omega(r)\ge 0$  if
		$$\|f\|_{\mathcal{B}(\omega)}:=|a_0|+\sup_{z\in \mathbb{D}}\omega(|z|)|f^{\prime}(z)|< \infty.$$
		For the standart weight $\omega(r)=1-r^2$ we just write $\mathcal{B}$.
	\end{definition}

	We start with the question of finding the value of $R_{\mathcal{B} \rightarrow H^{\infty}}$. This problem has already been studied in \cite{KPS}, where the following lower and upper estimates are obtained: $0.624162 \ldots \geq R_{\mathcal{B} \rightarrow H^{\infty}} \geq R$, where $R \approx 0.55356 \ldots$ is the solution of the equation $1-r+r \log (1-r)=0$. We noticed that the proof of the theorem can be slightly improved and we get a lower estimate about $0.01$ better. Note that in \cite{AHal} V. Allu and H. Halder considered this problem for functions defined on simply connected domains. As a consequence, they obtained a slightly less sharp lower estimate for the Bohr radius from the Bloch space to $H^{\infty}: R_{\mathcal{B} \rightarrow H^{\infty}} \geq 0.546679$.

	\begin{theorem}\label{Th:1}
		$$R_{\mathcal{B}\rightarrow H^{\infty}}\ge 0.563777.$$
	\end{theorem}

	\vspace{2mm}

	Next, we consider the problem of finding the Bohr radius of the weighted Bloch  spaces. The following theorem immediately follows from the Cauchy inequality.

	\begin{theorem}\label{Th:2} 
		
		Let $\mathcal{B}(\omega)$ be Bloch space with the weight $\omega(r)\ge 0$.
		Then
		$$R_{\mathcal{B}(\omega)}\ge 1/\sqrt{2}.$$

	\end{theorem}

	We ask a natural question: for a given weight, is the Bohr inequality with $R=1/\sqrt{2}$ sharp? Using the method of E. Bombieri and J. Bourgain \cite{Bombb}, we give a criterion for the sharpness of such an inequality.
	
	\begin{theorem}\label{Th:3}
		
		The inequality
		$$\Big|\Big|\sum_{n\ge 0} |a_n|\Big(\dfrac{z}{\sqrt{2}}\Big)^n\Big|\Big|_{\mathcal{B}(\omega)} \le ||f||_{\mathcal{B}(\omega)} \eqno(1)$$
		is sharp, i.e. turns to equality for some nonzero $f\in \mathcal{B}(\omega)$, if and only if there exists $r_0\in\Big[\dfrac{1}{\sqrt{2}}, 1\Big]$, such that for all $r\in[0,1)$ the inequality  $$\dfrac{\omega(r)}{\omega(r_0)}\le \min{\Big(2-\dfrac{r}{r_0}; \dfrac{\sqrt{2}r_0+r}{\sqrt{2}r+r_0}\Big)}. \eqno(2) \\$$ 
		holds. Extreme functions up to multiplication by a constant have the following form:
		$$\int_{c}^{z}\dfrac{z/r_0-e^{i\phi}/\sqrt{2}}{1-e^{-i\phi}z/(\sqrt{2}r_0)}dz+C, \quad \phi \in [0,2\pi), \: c\in \mathbb{D}, \: C\in \mathbb{C}.$$
	\end{theorem} 
	
	\vspace{1mm}

	Using Theorem 3, we can find the examples of the weights for which the inequality (1) is sharp. We fix $r_0\in [1/\sqrt{2},1]$ and denote $\omega_1(r):=2-\dfrac{r}{r_0}, \: \omega_2(r):=\dfrac{\sqrt{2}r_0+r}{\sqrt{2}r+r_0}.$ Note that $\omega_2(r)\le \omega_1(r)$ for $r\le r_0$ and $\omega_1(r)\le \omega_2(r)$  for $r\ge r_0$. Let $\tilde{\omega}(r):=\dfrac{\omega(r)}{\omega(r_0)},$ then $\tilde{\omega}(r)\le h(r):=\min{(\omega_1(r), \omega_2(r))}$ and $\tilde{\omega}(r_0)=1.$ It is easy to see, that $h(r)$ is a decreasing and convex function for $r\in [0,r_0).$ So we can drawn a plot of the function $h(r):$
	
	\vspace{2mm}

	\begin{figure}[h!]
		\centerline{\includegraphics[width=0.45\linewidth, height=6.6cm]{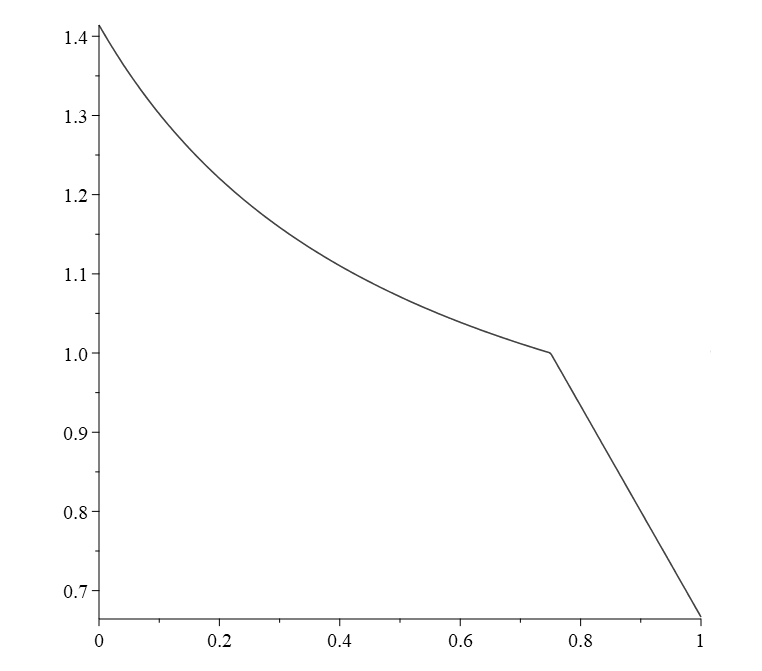}}  
	\end{figure}

	We see that there is no smooth function $\tilde{\omega}(r)$ with the necessary properties and $r_0<1$. That is why any smooth weight does not satisfy the condition (2). However, we can give examples of weights that satisfy condition (2):
	
	\vspace{4mm}{
		\textbf{Examples.}
		
		\vspace{1mm}
		1. $\omega(r)\equiv 1, \: r_0=1.$
		
		\vspace{3mm}
		Fix $r_0 \in [1/\sqrt{2}, 1)$ and $\alpha \ge 1.$
		
		\vspace{2mm}

		2. The weight 
		\begin{equation*}
			\omega(r) = 
			\begin{cases}
				1, & r\le r_0\\
				
				\Big(\dfrac{1-r}{1-r_0}\Big)^{\alpha}, & r\ge r_0
			\end{cases}
		\end{equation*}
		
		satisfies the condition (2).
		
		\vspace{2mm}
		3. The weight 
		
		$$\omega(r)=\Big(1-\Big|\dfrac{r-r_0}{1-r_0r}\Big|\Big)^{\alpha}$$
		
		also satisfies the condition (2).
		
		\vspace{2mm}
		Next we get the upper bound for the Bohr radius in classical Bloch space, using a function from the work of 	F. Avkhadiev and I. Kayumov \cite{Avkay}.

		\begin{theorem}\label{Th:4}
			Let $\mathcal{B}$ be the classical Bloch space with the weight $\omega(r)=1-r^2$. Then 
			$$ 0.707...\le R_\mathcal{B}\le 0.769...$$
		\end{theorem} 
		
		\vspace{1mm}
		
		\begin{problem}
			Let $X$ be Banach space of analytic functions in the unit disc, $f(z)=\sum_{n\ge 0}a_nz^n$. We need to find 
			$$m_X(r):=\sup_{f\in X}\dfrac{\Big\|\sum_{n\ge 0}|a_n|(rz)^n\Big\|_X}{\|f\|_X}.$$
		\end{problem}

		\vspace{3mm}
		In 1962 E. Bombieri showed (\cite{Bomb} or  \cite[p. 170--178]{GMR}) that for  $X=H^{\infty}$
		$$m_{\infty}(r):=m_{H^{\infty}}(r)=\dfrac{3-\sqrt{8(1-r^2)}}{r}, \quad 1/3\le r\le 1/\sqrt{2}.$$
		This theorem was proved by I. Kayumov and S. Ponnusamy in a simpler way in the work \cite{Kayump}.
		
		It immediately follows from the Cauchy inequality that $m_{{\infty}}(r)\le \dfrac{1}{\sqrt{1-r^2}}:$
		$$\sum_{n\ge 0}|a_n|r^n\le \Big({{\sum_{n\ge 0}|a_n|^2}}\Big)^{1/2}\Big({\sum_{n\ge 0}r^{2n}}\Big)^{1/2}=\dfrac{\|f\|_2}{\sqrt{1-r^2}}\le  \dfrac{\|f\|_\infty}{\sqrt{1-r^2}}.$$
		In the 2004, E. Bombieri and J. Bourgain proved \cite{Bombb} the following two profound theorems related to the inequality above.
		
	\newtheorem*{B}{Theorem B}\label{Th:B}
\begin{B}
			For all $r\in(0;1)$ except $r=1/\sqrt{2}$
			$$m_{{\infty}}(r)<\dfrac{1}{\sqrt{1-r^2}}, \:\:\:\:\: m_{{\infty}}(1/\sqrt{2})=\sqrt{2}.$$
			
		\end{B}
		
	\newtheorem*{C}{Theorem C}\label{Th:C}
	\begin{C}
			For all $\varepsilon>0$ there exists $C(\varepsilon)>0$, such that as $r \rightarrow 1$ the lower bound
			$$m_{{\infty}}(r)\ge \dfrac{1}{\sqrt{1-r^2}}-C(\varepsilon)\Big(\log\dfrac{1}{1-r}\Big)^{\varepsilon+3/2}$$ holds
		\end{C}
		
		Theorem C implies that $m_{{\infty}}(r)\rightarrow \infty$ when $ r\rightarrow 1$.
		It is clear from the proof of Theorem 1 that $m_\mathcal{B}(r)\le \dfrac{r}{\sqrt{1-r^2}}$. Using the construction already discussed in the proof of Theorem 2, we prove an analogue of Theorem B for the Bloch space:

		\begin{theorem}
			$$m_\mathcal{B}(R)<\dfrac{R}{\sqrt{1-R^2}}, \quad 0<R<1.$$
		\end{theorem}

		\vspace{4mm}
		\section{Proofs}
		
		\vspace{2mm}
		\begin{proof}[\textbf{Proof of Theorem 1}]
			If $\|f\|_\mathcal{B}\le1$, then $|f^{\prime}(z)|^2\le \dfrac{1}{(1-|z|^2)^2}, \: z\in \mathbb{D}$.
			
			\vspace{2mm}
			We integrate this inequality over the circle $|z|=r$ and obtain from Parceval formula $$\sum_{n\ge 1}n^2|a_n|^2r^{2n}\le \dfrac{r^2}{(1-r^2)^2},\quad0<r<1.$$
			Hence
			$$\sum_{n\ge 1}n^2|a_n|^2r^{n}\le \dfrac{r}{(1-r)^2}.$$
			We may now divide this inequality by r and integrate with respect to $r$ with limit from $0$ to $r$:
			$$\sum_{n\ge 1}n|a_n|^2r^n\le \dfrac{r}{1-r}.$$
			Replace $r$ with $r^s$:
			$$\sum_{n\ge 1}n|a_n|^2r^{ns}\le \dfrac{r^s}{1-r^s}, \quad 0<r, s<1.$$
			Thus,
			$$\sum_{n\ge 1}|a_n|r^n  \le\Big({\sum_{n\ge 1}\dfrac{r^{2ns}}{n}}\Big)^{1/2}\Big({\sum_{n\ge 1}n|a_n|^2r^{2n(1-s)}}\Big)^{1/2}\le\Big({\dfrac{r^{2(1-s)}\log(1-r^{2s})}{r^{2(1-s)}-1}}\Big)^{1/2}.$$
			Consequently, $\sum_{n\ge 1}|a_n|r^n\le 1$ with $r\le R$, where $R$ is solution of the equation 
			$$\log(1-r^{2s})=\dfrac{r^{2(1-s)}-1}{r^{2(1-s)}}.$$
			Using Wolfram Mathematica, we deduce that the best possible $r$, which can be achieved in this way is obtained at $s=0.333771$ and is equal to $0.563777.$
		\end{proof}
		
		\vspace{3mm}
		Unfortunately, there is no extreme function for $R=0.563777$, since if such a function $f(z)=\sum_{n\ge 1}a_nz^n$ existed, then
		$$|a_n|=C\dfrac{R^{n(2s-1)}}{n}, \quad n\ge 1$$
		Then the function $f(z)$ is analytic in $\mathbb{D}$ only for $s\ge 1/2$.
		
		\vspace{7mm}
		\begin{proof}[\textbf{Proof of Theorem 2}]
			In what follows we will solve an equivalent problem for weighted norms of functions in place of their derivatives. Note that $(f(Rz))' = Rf'(Rz)$. Thus, to find the Bohr radius in the Bloch space, it is sufficient to solve the following problem: find the maximal $R$ at which	
			$$	\sup_{{\substack{ r\in (0;1) }} }\omega(r)\|f(re^{i\theta})\|_{L^{\infty}[0,2\pi)}\le 1 \Longrightarrow  
			\\\sup_{r\in (0;1)}\omega(r)R\sum_{n\ge0}|a_n|(Rr)^n\le1,$$
			
			where $f(z)=\sum_{n\ge 0}a_nz^n$.
			
			\vspace{4mm}
			Let us show that $ R_{\mathcal{B}(\omega)}\ge \dfrac{1}{\sqrt{2}}:$
			\begin{multline*}
				\omega(r)R\sum_{n\ge 0}|a_n|(Rr)^n\le 
				\omega(r)R\Big({\sum_{n\ge 0}|a_n|^2r^{2n}}\Big)^{1/2}\Big({\sum_{n\ge 0}R^{2n}}\Big)^{1/2}= \label{eqn:line-3} \hspace{52mm} (4) \\ = \omega(r)\|f(re^{i\theta})\|_{L^2[0,2\pi]}\dfrac{R}{\sqrt{{1-R^2}}} \le\omega(r)\|f(re^{i\theta})\|_{L^{\infty}[0,2\pi]}\dfrac{R}{\sqrt{{1-R^2}}}. \notag
			\end{multline*}
			
			The inequality	$\dfrac{R}{\sqrt{{1-R^2}}}\le 1$ is equivalent to $R\le \dfrac{1}{\sqrt{2}}$, hence  $ R_{\mathcal{B}(\omega)}\ge \dfrac{1}{\sqrt{2}}.$
		\end{proof}
		
		\vspace{10mm}
		
		\begin{proof}[\textbf{Proof of Theorem 3}]
			As before, instead of derivatives of functions, we consider the functions themselves. Thus, we need to show that a nonzero function $f(z)=\sum_{n\ge 0}a_nz^n$ such that
			$$\sup_{r\in (0;1)}\dfrac{\omega(r)}{\sqrt{2}}\sum_{n\ge 0} |a_n|\Big(\dfrac{r}{\sqrt{2}}\Big)^n=\sup_{{\substack{ r\in (0;1) }} }\omega(r)\|f(re^{i\theta})\|_{L^{\infty}[0,2\pi)}<\infty. \hspace{3mm}\eqno(5)$$
			exist if and only if the weight $\omega$ satisfies (2).
			
			\vspace{2mm}
			\textit{Necessity.} Assume that equality (5) holds for some $f$. Let us show that there exists $r_0\in [0,1]$ with the following properties:

			\vspace{3mm}
			1. $r_0$ is the point where the supremum of the left-side in (5) is attained:
			$$\omega(r)\sum_{n\ge 0} |a_n|\Big(\dfrac{r}{\sqrt{2}}\Big)^n\le \omega(r_0)\sum_{n\ge 0} |a_n|\Big(\dfrac{r_0}{\sqrt{2}}\Big)^n, \quad0\le r<1.$$
			
			2. The norm of $f$ is also attained at $r_0:$
			$$\omega(r)\|f(re^{i\theta})\|_{L^{\infty}[0,2\pi)}\le  \omega(r_0)\|f(r_0e^{i\theta})\|_{L^{\infty}[0,2\pi)}, \quad 0 \le r<1.$$

			3. The Cauchy inequality in the chain (4) turns into equality which is equivalent to  $$|a_n|r_0^n=C(r_0)\Big(\dfrac{1}{\sqrt{2}}\Big)^n,\quad n\ge 0,$$ where $C(r_0)$ is a constant depending only on $r_0$. 
			
			4. The last inequality in the chain (4) is equality, that is
			$$\|f(r_0e^{i\theta})\|_{L^{2}[0,2\pi)}=\|f(r_0e^{i\theta})\|_{L^{\infty}[0,2\pi)}.$$
			Define $r_0$ by the condition 1. We substitute $r=r_0$ and $R=1/\sqrt{2}$  in (4) and use the condition 1 and equality (5):
			\begin{multline*} \omega(r_0)\|f(r_0e^{i\theta})\|_{L^{\infty}[0,2\pi)}\ge \dfrac{\omega(r_0)}{\sqrt{2}}\sum_{n\ge 0} |a_n|\Big(\dfrac{r_0}{\sqrt{2}}\Big)^n \ge \\ \ge
				\sup_{r\in (0;1)}\dfrac{\omega(r)}{\sqrt{2}}\sum_{n\ge 0} |a_n|\Big(\dfrac{r}{\sqrt{2}}\Big)^n=	\sup_{r\in (0;1)}\omega(r)\|f(re^{i\theta})\|_{L^{\infty}[0,2\pi)}.	\end{multline*}
			Thus, condition 2 is also satisfied. Conditions 3 and 4 following from the fact that all inequalities in (4) are equalities for $r=r_0$ and $R=1/\sqrt{2}.$ 
			
			It follows from the condition 3 that the radius of convergence for $f$ is $\sqrt{2}r_0$, where $r_0\ge 1/\sqrt{2}.$
			It is known that the equality $\|f(r_0e^{i\theta})\|_{L^{2}[0,2\pi)}=\|f(r_0e^{i\theta})\|_{L^{\infty}[0,2\pi)}$ holds only for functions $$f(z)=K(r_0)B\Big(\dfrac{z}{r_0}\Big),$$ where the function $B(z)$ is a finite Blaschke product (see, for example, \cite{GMR}). Without loss of generality, we can assume that $|K(r_0)|=1$.

			We find $C(r_0)$. Let $$B\Big(\dfrac{z}{r_0}\Big)=\prod_{i=1}^{K}\Big(\dfrac{z/r_0-\alpha_i}{1-\bar{\alpha_i}z/r_0}\Big)^{d_i}, \hspace{7mm}\sum_{i=1}^{K}d_i=d, \hspace{5mm} \alpha_i\in \mathbb{D} \:(1\le i\le K).$$
			
			Denote by $dA(z)$ the standard Lebesgue measure on the plane. From the formula for the area of the image of the finite Blaschke product, we can calculate its degree:
			\begin{eqnarray*}d=
				\dfrac{1}{\pi}\int_{\mathbb{D}}\big|B^{\prime}(z)\big|^2 dA(z)=\dfrac{1}{\pi}\int_{\mathbb{D}}\big|(f(r_0z))^{\prime}\big|^2 dA(z) =\\ 
				=\sum_{n\ge 1}n|a_n|^2r_0^{2n}=\sum_{n\ge 1}\dfrac{n}{2^n}C^2(r_0)=2C^2(r_0).
			\end{eqnarray*} 
			Therefore $C(r_0)=\sqrt{d/2}.$
			Note that
			$$\sqrt{d/2}=C(r_0)=|a_0|=|f(0)|=\Big|\prod_{i=1}^{K}\alpha_i^{d_i}\Big|<1, \eqno(6)$$
			hence $d=1$.
			If $d=1$, then from (6), $|\alpha_1|=1/\sqrt{2}$ and$$f(z)=\dfrac{z/r_0-e^{i\phi}/\sqrt{2}}{1-e^{-i\phi}z/(\sqrt{2}r_0)}.$$
			Without loss of generality, we consider $\phi=0$. In this case, the supremum $|f(re^{i\theta})|$ is reached when $\theta=\pi$,
			$$\|f(rz)\|_{L^\infty[0,2\pi)} = |f(-r)|=\dfrac{r/r_0+1/\sqrt{2}}{1+r/(\sqrt{2}r_0)}$$ 
			and
			$$\sum_{n\ge 0} |a_n|\Big(\dfrac{r}{\sqrt{2}}\Big)^n=\dfrac{1}{\sqrt{2}}\dfrac{1}{1-r/(2r_0)}.$$
			
			Now it follows from the condition 1 that 
			$$\dfrac{\omega(r)}{\omega(r_0)}\le 2-\dfrac{r}{r_0},$$
			while condition 2 implies
			$$\dfrac{\omega(r)}{\omega(r_0)}\le \dfrac{\sqrt{2}r_0+r}{\sqrt{2}r+r_0}.$$
			\vspace{0.00000000000000000003mm}
			
			We showed that the weight $\omega$ satisfies (2).
			
			\vspace*{5mm}
			\textit{Sufficiency}.
			Define $f$ by 
			$$f(z)=\dfrac{z/r_0-e^{i\phi}/\sqrt{2}}{1-e^{-i\phi}z/(\sqrt{2}r_0)}.
			$$
			Then $f$ automatically satisfies condition 3 and 4, while condition (2) on the weight implies that $f$ also satisfies condition 1 and 2. Hence for $f$ condition 1 holds with both suprema attained at $r=r_0.$\end{proof}
		
		\vspace{6mm}
		\begin{proof}[\textbf{Proof of Theorem 4}] We consider the function $f(z)=\sum_{n\ge 0}a_nz^n$ from \cite{Avkay}:
			$$f(z)=\dfrac{3\sqrt{3}}{2}\dfrac{1-a^2}{(1-az)^2}\Big(\dfrac{z-a}{1-az}\Big), \quad 0<a<1/\sqrt{3}.$$
			It is not difficult to show that $\sup_{z\in\mathbb{D}}(1-|z|^2)|f(z)|=1$. Obviously, $a_0 = -\frac{3\sqrt{3}a(1-a^2)}{2}<0$. Let us show that if $0<a<1/\sqrt{3}$, then $a_n>0, \:n\ge 1$. Denote $\tilde{a}=\dfrac{a^2}{1-a^2}$ and let
			$$g(w)=\sum_{n\ge 0}b_nw^n=\dfrac{1}{(1-w)^2}\Big(\dfrac{w}{1-w}-\tilde{a}\Big)=\sum_{n\ge 0}(n+1)(n/2-\tilde{a})w^n.$$
			Since $f(z)=C(a)g(az), \:C(a)>0,$  the condition $a_n>0$ is equivalent to $b_n>0.$ But $b_n=(n+1)(n/2-\tilde{a})>0$ for all $n\ge 1$ if and only if $0<a<1/\sqrt{3}.$
			
			\vspace{1mm}
			Thus,
			$$\sum_{n\ge 0}|a_n|z^n=  \dfrac{3\sqrt{3}(1-a^2)}{2}\Big(\dfrac{z-a}{(1-az)^3}+2a\Big)$$ and
			$$R(1-r^2)\sum_{n\ge 0}|a_n|(Rr)^n=R(1-r^2)\dfrac{3\sqrt{3}(1-a^2)}{2}\Big(\dfrac{Rr-a}{(1-aRr)^3}+2a\Big).$$
			Using Wolfram Mathematica again, we find that the last expression is greater than one, for example, at $a=0.35$ and $R=0.769$, so $R_\mathcal{B}\le 0.769$.
		\end{proof}
		
		\vspace{9mm}
		\begin{proof}[\textbf{Proof of Theorem 5}]
			Suppose that the expression $(1-r^2)\sum_{n\ge 0} |a_n|\Big(\dfrac{r}{\sqrt{2}}\Big)^n$ reaches the supremum at the point $r_0$. We replace all the  inequalities in the chain (4) with equalities for $r=r_0$. Then, as in the proof of Theorem 3, it will follows from the sharpness condition of the Cauchy's inequality that
			$$|a_n|r_0^n=C(r_0)R^n,\quad n\ge 0,$$ so
			$r_0\ge R$.  It follows from the equality of norms that $f(z)=B\Big(\dfrac{z}{r_0}\Big)$, where $B(z)$ is finite Blaschke product that is
			$$f(z)=B\Big(\dfrac{z}{r_0}\Big)=\prod_{i=1}^{K}\Big(\dfrac{z/r_0-\alpha_i}{1-\bar{\alpha_i}z/r_0}\Big)^{d_i}, \:\:\: \sum_{i=1}^{K}d_i=d.$$ 
			It follows from the formula for the area of the image of the finite Blaschke product  that
			$$C(r_0)=\sqrt{d} \Big(\dfrac{1-R^2}{R}\Big).$$
			
			Since we have replaced all the inequalities with equalities, then
			
			$$(1-r^2)|f(r)|\le (1-r_0^2)|f(r_0)|=1-r_0^2, \quad 0 \le r<1. $$
			We substitute $r=0$ into the left part of the inequality above:
			$$|f(0)|=C(r_0)=\sqrt{d}\Big(\dfrac{1-R^2}{R}\Big)\le 1-r_0^2\le1-R^2,$$
			which is impossible for all $R\in(0;1)$.
		\end{proof}
		
	\textbf{Acknowledgments.}
			I would like to thank  Professor A.D. Baranov for formulation of the problem and for valuable discussions, Professor I.R. Kayumov for useful comments and M.S. Ivanova for help with examples of weights.
			
			The work is supported by Ministry of Science and Higher Education of the Russian Federation (agreement No 075-15-2021-602) and by Theoretical Physics and Mathematics Advancement Foundation "BASIS".}
		
		%
		%

	\end{document}